\documentclass[12pt]{article}
\usepackage{fullpage,graphicx,amsmath,amsfonts}
\usepackage[small,bf]{caption}
\usepackage{authblk}
\usepackage{url}
\usepackage{booktabs}
\usepackage{hyperref}

\hypersetup{
    colorlinks=true,     
    linkcolor=blue,      
    citecolor=red,       
    urlcolor=magenta     
}

\newcommand{\reals}{{\mbox{\bf R}}}

\newcommand{\eg}{{\it e.g.}}
\newcommand{\ie}{{\it i.e.}}

\newcommand{\BEAS}{\begin{eqnarray*}}
\newcommand{\EEAS}{\end{eqnarray*}}
\newcommand{\BEA}{\begin{eqnarray}}
\newcommand{\EEA}{\end{eqnarray}}
\newcommand{\BEQ}{\begin{equation}}
\newcommand{\EEQ}{\end{equation}}
\newcommand{\BIT}{\begin{itemize}}
\newcommand{\EIT}{\end{itemize}}
\newcommand{\BNUM}{\begin{enumerate}}
\newcommand{\ENUM}{\end{enumerate}}

\newcommand{\BEQN}{\begin{equation*}}
\newcommand{\EEQN}{\end{equation*}}

\newcounter{algorithmctr}
\renewcommand{\thealgorithmctr}{\arabic{algorithmctr}}
   {\mbox{}\\*[\parskip]\begin{minipage}{\linewidth}%
       \refstepcounter{algorithmctr}\begin{list}{}{%
       \setlength{\rightmargin}{0\linewidth}%
       \setlength{\leftmargin}{.05\linewidth}}%
       \rmfamily\small
       \item[]{\setlength{\parskip}{0ex}\hrulefill\par%
        \nopagebreak{\bfseries\textsf{Algorithm \thealgorithmctr~}}}}%
   {{\setlength{\parskip}{-1ex}\nopagebreak\par\hrulefill\\*[2ex]\par}%
   \end{list}\end{minipage}}

\title{GPU-Accelerated Presolving for Linear Programming}
\author[]{Daniel Cederberg \qquad \qquad Stephen Boyd \\ Stanford University}

\begin{document}
\maketitle

\begin{abstract}
Recent research has focused on developing GPU-accelerated algorithms for solving
linear programs (LPs), with results that are nothing short of extraordinary. But
a complete solver pipeline consists of more than the core algorithm, and the
elephant in the room is that stages such as presolve have remained on the CPU,
where they have become an increasingly large bottleneck. Presolve has stayed on
the CPU because of the longstanding view, held by academic and industrial
developers alike, that it is inherently sequential and irregular, and therefore
difficult to parallelize efficiently. In this paper we challenge this view and
show that presolving can benefit substantially from GPU acceleration. We
describe simple design principles that avoid severe load imbalance and expose
parallelism in reductions that appear inherently sequential. We implement these
principles in cuPSLP, a GPU-accelerated version of PSLP, a CPU-based presolver
that is itself several times faster than a state-of-the-art commercial
presolver. In our experiments, cuPSLP reduces the shifted geometric mean
presolve time of PSLP by a factor of 11 on the Mittelmann LP benchmark set and
42 on the GAMS large-scale LP benchmark set, with similar reduction quality.
\end{abstract}

\newpage
\tableofcontents
\newpage
\section{Introduction}
Linear programming (LP) is among the most widely used tools in mathematical
optimization, with applications spanning a broad range of fields. An important
component of modern LP solvers is the \emph{presolver} \cite{Brearley1975,
Andersen1995a, Achterberg2019}, which fixes variables and removes redundant
constraints to shrink the problem before a \emph{core algorithm} (such as an
interior-point method) solves the reduced problem.

Recently developed solvers have focused on accelerating the core algorithm by
leveraging modern parallel hardware. Examples include PDLP \cite{Applegate2025},
a multithreaded CPU implementation of the primal-dual hybrid gradient method
\cite{Zhu2008, Chambolle2011}, and its GPU descendants in the cuPDLP family
\cite{Lu2025, Lu2023, Lu2025-cuPDLPx}. While these solvers exploit massive
parallelism in their core algorithm, their presolve step, if any, does not. For
example, cuPDLPx \cite{Lu2025-cuPDLPx}, NVIDIA's GPU-accelerated library cuOpt
\cite{CuOpt2025}, and HPR-LP \cite{Chen2025} all integrate the presolver PSLP
\cite{Cederberg2026}, which runs (mostly) single-threaded on the CPU before the
reduced problem is transferred to the GPU only for the solve.

Presolving has traditionally accounted for a small fraction of the total solve
time. But as GPU solvers have made the core algorithm extraordinarily fast, the
cost of a sequential CPU presolve can no longer hide behind it. Even PSLP, a
lightweight presolver engineered for speed, accounts on average for 28\% of
end-to-end solve time for cuPDLPx on a subset of the larger LP relaxations of
the MIPLIB 2017 benchmark set. For a more complex presolver such as Gurobi's
\cite{Achterberg2019}, this proportion rises to 66\% \cite[Table
4]{Cederberg2026}. One might respond by simply skipping presolve, but presolve
is often beneficial and sometimes essential to solving a problem at all, and
thus a critical part of any complete solver pipeline. As the core algorithm
keeps getting faster, however, a presolve that is itself not accelerated will
consume an ever-larger share of the total solve time and may come to dominate
it.

The prevailing view, held by some of the world's best optimization developers,
is that presolve is an inherently sequential task, with little room for
parallelism and in particular for GPU acceleration \cite{Rothberg2025}. (The
scarcity of work on parallel presolving, discussed below, is further evidence
that this view is widely held.) We too shared this view in our previous work
\cite{Cederberg2026}, for the following reasons. To realize the potential of a
GPU, a computation must expose large amounts of fine-grained parallelism and
exhibit sufficient regularity of execution paths and memory access patterns.
Presolving does neither well: it proceeds through long sequences of
transformations that appear inherently sequential, and it is routinely applied
to problems whose rows and columns vary by orders of magnitude in length,
leading to severe load imbalance if parallelized naively. (We discuss this in
more detail in \S \ref{sec:load-imbalance}.) These difficulties are typical of
sparse linear algebra, which some presolve operations resemble. Even highly
tuned kernels for sparse linear algebra reach only a small fraction of a GPU's
peak throughput \cite{Mpakos2023}, and presolving is even more irregular.
Largely for these reasons, presolving has remained on the CPU and come to be
regarded as an inherently CPU-bound task.

The central finding of this paper is that, contrary to this prevailing view,
\emph{presolving can greatly benefit from GPU acceleration}. We present simple
design strategies that avoid severe load imbalance and expose parallelism in
reductions that appear inherently sequential, and we implement them in cuPSLP, a
GPU-accelerated version of the CPU-based presolver PSLP \cite{Cederberg2026}. In
our experiments, cuPSLP running on an A100 presolves more than one order of
magnitude faster than PSLP running on a Xeon CPU (see \S \ref{sec:experiments}
for the exact hardware specifications). This matters because presolving is one
of the last sequential, CPU-side stages of an otherwise GPU-resident solver
pipeline. Accelerating it removes a bottleneck that only grows as the core
algorithm gets faster. Moreover, the several solvers that already rely on PSLP
can adopt cuPSLP with little effort and benefit immediately.

\paragraph{Additional related work.}
The idea of presolving an optimization problem before passing it to a solver is
half a century old. Early work \cite{Brearley1975, Williams1983, Bradley1983,
Tomlin1986} focused primarily on linear programming with the simplex method as
the core algorithm. The development of interior-point methods in the 1990s
renewed interest in presolving \cite{Andersen1995a, Andersen1995b, Gondzio97,
Meszaros2003, Gould2004}. Most recently, we revisited its role in the context of
GPU-accelerated first-order solvers \cite{Cederberg2026}.

A few efforts have aimed to parallelize presolving itself. PaPILO
\cite{Gleixner2023} is a multithreaded CPU presolver for mixed-integer linear
programs that uses a transaction-based design, allowing independent reductions
to proceed concurrently. The presolver of Kempke et al.\ \cite{Kempke2026} runs
across distributed-memory compute nodes and is designed specifically for LPs
with arrowhead structure. On GPUs, Sofranac et al.\ \cite{Sofranac2022}
parallelize bound tightening, which is a single presolve technique, rather than
a complete presolver. The relatively scarce literature on parallel presolving is
consistent with the prevailing view that presolving is inherently sequential,
and that parallelism is difficult to exploit.

\paragraph{Outline.} The remainder of the paper begins with a brief overview of
presolve and the GPU execution model in \S \ref{sec:background}. In
\S \ref{sec:challenges}, we discuss important design considerations for a
GPU-accelerated presolver, and give a brief overview of how we implement some of
the GPU-accelerated exploration kernels within cuPSLP. We present numerical
results in \S \ref{sec:experiments}, followed by a discussion in
\S \ref{sec:discussion}.

\section{Background}
\label{sec:background}

\subsection{Presolve}
\label{sec:bg-presolve}
We consider an LP of the form
\begin{equation}\label{eq:lp}
\begin{array}{ll}
\mbox{minimize} & c^T x \\
\mbox{subject to} & \underline{b} \leq Ax \leq \overline{b} \\
& \underline{x} \leq x \leq \overline{x},
\end{array}
\end{equation}
with variable $x \in \reals^n$. The problem data are the cost vector $c \in
\reals^n$, the left- and right-hand sides $\underline{b} \in \reals^m$ and
$\overline{b} \in \reals^m$, the constraint matrix $A \in \reals^{m \times n}$,
and the variable bounds $\underline{x} \in \reals^n$ and $\overline{x} \in
\reals^n$. Infinite values are allowed in both the variable bounds and the left- and
right-hand sides, allowing for unbounded variables and one-sided inequality
constraints. Equality constraints can be expressed by setting the corresponding
components of $\underline{b}$ and $\overline{b}$ equal.

We now describe how a presolver simplifies \eqref{eq:lp}, using the abstractions
of \emph{reductions} and \emph{explorers}, and how it reconstructs a solution to
the original problem using a \emph{postsolver}. For more background on presolve,
we refer to \cite{Andersen1995a, Achterberg2019, Cederberg2026}.

\paragraph{Reductions.}
A presolver simplifies \eqref{eq:lp} by applying a sequence of reductions, each
of which consists of a \emph{presolve transformation} and a \emph{postsolve
transformation}. The presolve transformation takes as input one problem and
outputs an equivalent problem that is in some sense smaller or simpler. The
postsolve transformation takes as input a primal-dual solution to the reduced
problem and maps it to a primal-dual solution of the input problem. Each
reduction is simple and cheap to apply, and for large LPs several hundred
thousand reductions may be applied in sequence. Examples of reductions are
substituting a variable, fixing a variable to some value, or removing a
constraint.

\paragraph{Explorers.}
To identify opportunities for invoking reductions, a presolver applies a set of
explorers. Each explorer implements logical rules based on either feasibility or
optimality: \emph{primal explorers} identify reductions that preserve the primal
feasible set, while \emph{dual explorers} identify reductions that preserve at
least one optimal solution. An example is the doubleton-row explorer, which
detects equality constraints with two nonzeros and invokes a reduction that uses
the constraint to eliminate one of the two variables. Further examples of
explorers are described in \cite{Achterberg2019}.

Designing an explorer involves a tradeoff between the cost of applying it and
the size of the reduction it can achieve. Many explorers have been proposed,
varying widely in cost, but in previous work \cite{Cederberg2026} we showed that
a carefully engineered collection of relatively simple ones recovers most of the
reduction (on average about 90\%) obtained by a more complex presolver such as
Gurobi's \cite{Achterberg2019}.

\paragraph{Postsolver.}
Each time a reduction is applied during presolve, the data required for its
postsolve transformation is pushed onto a \emph{postsolve stack}. Given a
primal-dual solution to the reduced problem, the postsolver pops this stack one
reduction at a time to reconstruct a primal-dual solution to the original
problem. Each reduction's postsolve transformation maps the current solution
back to a solution of the problem as it stood just before that reduction was
applied.

The amount of data that must be stored varies greatly between reductions. For
example, fixing a variable requires storing only its index and value, while
substituting a variable using an equality constraint requires storing the entire
constraint itself (so that the value of the eliminated variable can be
reconstructed from the optimal values of the remaining variables).

\subsection{The GPU execution model}
\label{sec:bg-gpu}
We now describe the GPU execution model, covering only the aspects we rely on
later in the paper. There is considerably more to the execution model than we
present here, but we have found these basics sufficient to achieve substantial
speedups over CPU-based presolve.

\paragraph{Warps and thread blocks.}
A function that runs on the GPU is called a \emph{kernel} and is executed by a
large number of threads in parallel that are grouped into \emph{warps} of $32$
(on NVIDIA GPUs). A warp essentially executes a single instruction at a time
across all its threads, so it proceeds at the pace of its slowest thread,
and threads that finish early sit idle until the rest catch up. When threads
within a warp follow different execution paths (\ie, branches of if-statements),
the warp executes each branch in turn with the inactive threads masked off, a
phenomenon known as \emph{warp divergence}. Threads are organized by the
programmer into \emph{thread blocks}, which the hardware in turn partitions into
warps. Threads within a block share a fast on-chip memory and can synchronize
with one another. Thread creation and warp scheduling are performed entirely by
the GPU hardware, and the programmer controls them only indirectly by specifying
how many threads to launch and how to distribute work among warps and blocks.

\paragraph{CPU versus GPU threads.}
A single GPU thread is much slower than a single CPU thread (by which we mean
one sequential instruction stream on one CPU core). One reason is that a CPU
core hides memory latency within one thread, using, \eg, caches and prefetching,
while a GPU hides it by switching among many resident warps. With only one warp
resident, there is nothing to switch to, so the full latency of every memory
access is exposed. For example, summing an array of $2.4$ million doubles (the
length of the longest row of \texttt{Dagolpert}, a problem instance we discuss
below) takes 174 ms on a single A100 thread, versus 175 $\mu$s on a single Apple
M4 Pro core, making the lone GPU thread three orders of magnitude slower. (For
this experiment we use a plain C for-loop for the sum, average the time over
1000 runs, and flush the cache before each run. We compile both with
\texttt{-O3} and the CPU additionally with \texttt{-ffast-math} so that its
single thread can vectorize across the core's SIMD lanes, which on a GPU are
instead spread across a warp's threads.) The lesson is not that the GPU is slow
but that its performance comes from parallelism: work concentrated on a single
thread, as happens under load imbalance, leaves the rest of the device idle.

\paragraph{Atomics and warp reductions.}
When many threads update the same memory location concurrently, their updates
can interfere and corrupt the result. An \emph{atomic operation} performs a
read-modify-write that is guaranteed to complete without interruption from other
threads, and is the standard way to make such concurrent updates safe. Atomic
operations serialize the updates, making them slow if many threads access the
same location. A \emph{warp reduction} avoids this by having each of the 32
threads in a warp accumulate a partial result in registers, which are then
combined through warp shuffle instructions (each thread reads another lane's
register directly, with no memory traffic), so that only a single write per warp
reaches the shared location.

\section{Designing a GPU-accelerated presolver}
\label{sec:challenges}
In this section we discuss important design considerations for a GPU-accelerated
presolver, and give a brief overview of how we implement some of the
GPU-accelerated exploration kernels within cuPSLP.

\subsection{Load imbalance}
\label{sec:load-imbalance}
\begin{figure}[ht]
\centering
\includegraphics[width=0.7\textwidth]{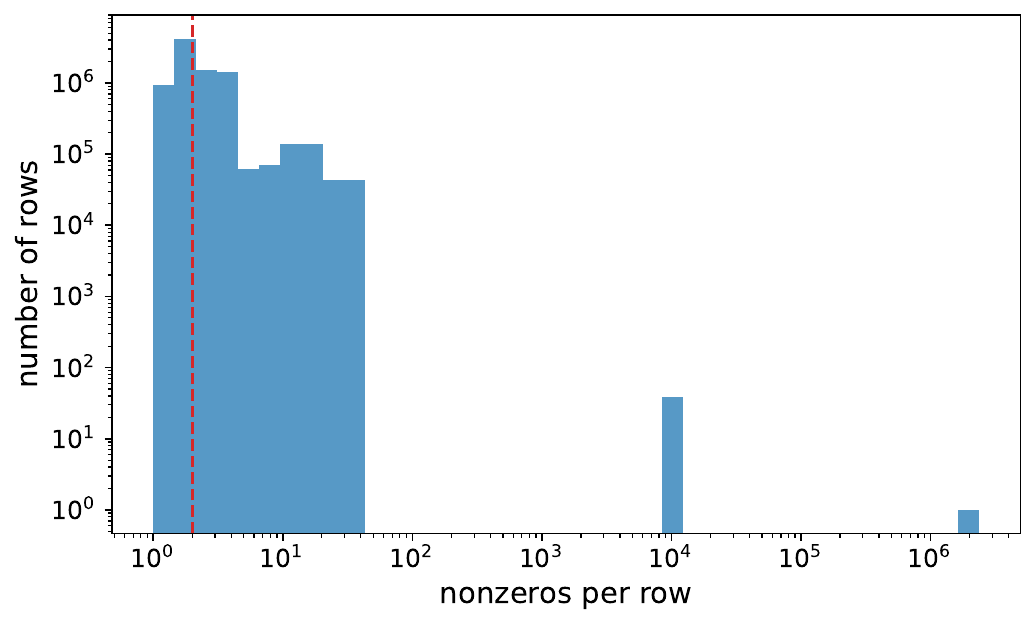}
\caption{Row-length distribution of the GAMS instance \texttt{Dagolpert}. The
red dashed line marks the median row length of $2$ nonzeros.}
\label{fig:rowdist}
\end{figure}

Load imbalance is a classic challenge for sparse linear algebra on GPUs
\cite{Bell2009, Greathouse2014}. Presolving inherits this challenge, since many
of the operations performed by a presolver are sparse linear algebra operations.
(For example, bound tightening, which we discuss in more detail below, is
essentially a sparse matrix-vector product.) The imbalance arises because the
work done by an algorithm on a row or column is typically proportional to its
length, and these lengths can vary by orders of magnitude. For example, Figure
\ref{fig:rowdist} shows the row-length distribution for the instance
\texttt{Dagolpert} from the GAMS large-scale LP benchmark set
\cite{GamsLp2026}. The median row has $2$ nonzeros, but the densest has nearly
$2.4$ million, and this longest row alone accounts for about 8\% of the matrix's
nonzeros. If one GPU thread is assigned to each row, the threads handling the
long rows take far longer than those handling the short rows, and the kernel
runs at the pace of its slowest thread. On a CPU this is far less of a concern:
a single-threaded implementation depends only on the total amount of work, not
its distribution, and even a multithreaded implementation suffers little, since
a long row that would stall a lightweight GPU thread is handled far more quickly
by a single powerful core.

Load imbalance can appear at many places in a presolver. Examples include
substituting a variable that occurs in a long row into the objective, storing
information about a long row on the postsolve stack, or scanning long rows to
tighten variable bounds. The simple yet powerful remedy is to assign more
threads to long rows, so that they are processed cooperatively rather than by a
single thread. We describe this below in the context of \emph{bound tightening},
but the same principle applies to other explorers and kernels as well.

\paragraph{Bound tightening.} Bound tightening, also referred to as
\emph{propagation}, is an explorer that uses the constraints together with
existing variable bounds to tighten other variable bounds (see, for example,
\cite{Achterberg2019, Sofranac2022}). For example, the constraint $x_1 + 2x_2
\leq 4$ together with the lower bounds $x_1 \geq 0$ and $x_2 \geq 0$ implies the
upper bounds $x_1 \leq 4$ and $x_2 \leq 2$, which in turn can be used to tighten
the bounds of variables in other constraints, and so on.

The implementation of bound tightening is based on computing the smallest and
largest value each constraint can take given the current variable bounds (the
so-called \emph{minimum} and \emph{maximum activities}), and then using these to
derive implied bounds on the variables appearing in the constraint. We carry
this out with a single fused kernel that assigns one unit of parallelism (\ie, a
thread, warp, or thread block, as discussed below) to each row: the unit first
computes the activities given the current bounds and then tightens the bounds of
the variables in the assigned row. We repeat this in rounds, stopping when a
round tightens no bounds, or after a small fixed number of rounds. Because the
bounds are updated in place (using atomic operations), a tightening made during
a round is seen immediately by the other rows of the same round, so the
iteration is effectively of Gauss--Seidel type \cite[\S 11.2]{Golub2013}.

\begin{figure}[ht]
\centering
\includegraphics[width=0.7\textwidth]{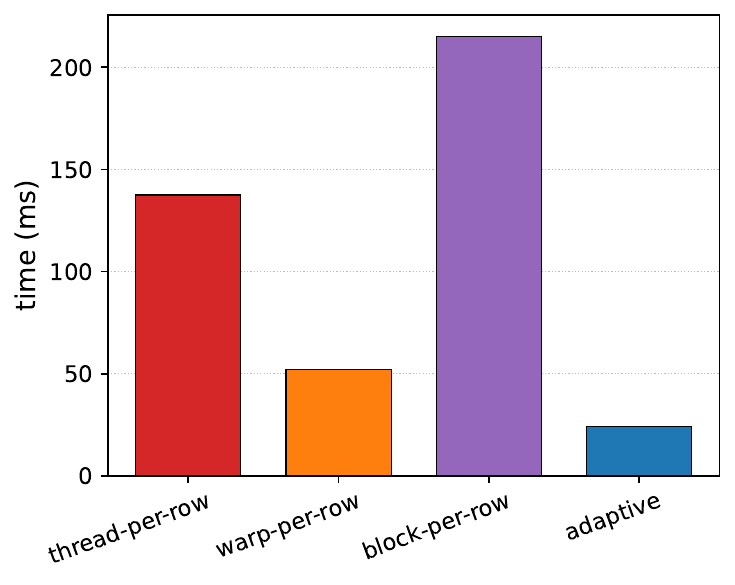}
\caption{Propagation time on \texttt{Dagolpert} with different load-balancing
strategies.}
\label{fig:loadbalance}
\end{figure}

\paragraph{Length-based load balancing.}
The cost of bound tightening on a row is proportional to its length, so the
distribution in Figure \ref{fig:rowdist} is also the distribution of work across
units of parallelism. Assigning one thread per row is inefficient: a thread that
draws a long row holds up the rest of its warp, which sits idle until it
finishes. Assigning a full warp or block to every row is also inefficient, since
short rows leave most of the threads with nothing to do. We instead use
length-based load balancing \cite{Greathouse2014}: we assign one thread to each
short row, one warp to each medium row, and one block to each long row, with the
boundaries chosen empirically. Before invoking the bound-tightening explorer, we
partition the active rows into these three groups, and then launch one kernel
per group. Figure \ref{fig:loadbalance} shows the effect on \texttt{Dagolpert}
on an A100 GPU, where the adaptive strategy is more than twice as fast as the
best static strategy. No single static strategy is best across the GAMS
large-scale LP benchmark set, yet the adaptive strategy is the fastest of the
four strategies in Figure \ref{fig:loadbalance} on every problem, with a median
speedup of $2.3\times$ over the best static strategy.

We have found that simple strategies like this are effective for many of the
explorers, and that they are crucial for eliminating the worst load imbalance
that causes bottlenecks.

\subsection{Exploration kernels}
We now give a brief overview of how some of the other explorers are parallelized
within cuPSLP.

\paragraph{Internal data structures.} To implement explorers efficiently, a
CPU-based presolver maintains and incrementally updates several internal
statistics about the problem, including the minimum and maximum activity for
each constraint, and the number of constraints that prevent each variable from
increasing or decreasing to $\pm\infty$ (the so-called \emph{variable locks}).
Moreover, some CPU-based presolvers such as PSLP are designed to eliminate a
variable by physically removing it from every row in which it appears, modifying
the sparse constraint matrix $A$ itself. In the GPU-accelerated setting, we
depart from this in two ways. Rather than physically removing an eliminated
variable from $A$, we mark it as inactive and ignore it when streaming through
$A$, which we found both simpler and more efficient. Moreover, instead of
maintaining the activities and locks incrementally, we recompute them on demand,
immediately before the explorers that rely on them. As a result, cuPSLP
maintains only the row and column sizes, along with flags marking which rows
and columns are inactive. Furthermore, since some of the
explorers require fast access to rows while others require fast access to
columns, we store both $A$ and its transpose $A^T$ in compressed sparse row
(CSR) format.

\paragraph{Parallel rows.} Two constraints $\underline{b}_i \leq a_i^T x \leq
\overline{b}_i$ and $\underline{b}_j \leq a_j^T x \leq \overline{b}_j$ are
\emph{parallel} if $a_i = \lambda a_j$ for some $\lambda \in \reals$. Parallel
constraints can be merged into a single one, after adjusting their left- and
right-hand sides accordingly \cite[\S 5.2]{Achterberg2019}. To search for
parallel constraints, we implement a two-level hash procedure that hashes each
row based on its sparsity pattern and coefficients \cite[\S 3.5]{Andersen1995a}.
Hash values for different rows can be computed independently, and we carry this
out using one warp per row. Rows that are parallel will have the same hash
value, but two rows with the same hash value are not necessarily parallel
because of the possibility of hash collisions. Therefore, we then explicitly
compare rows whose hashes collide. Within each group of rows sharing a hash,
every thread of a warp checks whether one or more of its assigned constraints
are parallel to a designated \emph{leader} constraint, which is taken to be the
first equality constraint in the group if one exists and the first inequality
constraint otherwise.

\paragraph{Parallel columns.}
Two variables $x_j$ and $x_k$ are \emph{parallel} if $a_j = \lambda a_k$ and
$c_j = \lambda c_k$ for some $\lambda \in \reals$, where $a_j$ and $a_k$ here
denote the $j$th and $k$th columns of $A$. Parallel columns can be merged into a
single one, after updating the bounds of the merged variable \cite[\S
3.6]{Andersen1995a}. We detect parallel columns by searching for parallel rows
of the augmented matrix $\begin{bmatrix} c & A^T \end{bmatrix}$.

\paragraph{Column singletons.}
A \emph{column singleton} is a variable $x_j$ that appears in only one
constraint. When $x_j$ is \emph{free} (meaning it has no finite bounds), or when
its bounds are implied by the single constraint it appears in, we can eliminate
$x_j$ from the problem \cite[\S 3.2]{Andersen1995a}. We process column
singletons using one thread per column. However, the substitution of $x_j$ into
the objective function can be expensive when the constraint $x_j$ appears in has
many nonzeros, so we defer the substitution when the length of the row exceeds a
threshold, recording the substitution as a \emph{task} to be processed by a
second kernel that assigns one thread block per task. This ensures that the
substitution of $x_j$ using the long row is processed cooperatively by many
threads.

\subsection{Our implementation}
We have implemented the above design in cuPSLP, a GPU-accelerated version of the
CPU-based presolver PSLP \cite{Cederberg2026} available at
\begin{center}
\url{https://github.com/dance858/PSLP}.
\end{center}
cuPSLP takes in an LP of the form \eqref{eq:lp}, and emits a reduced LP of the
same form. Given an (approximate) primal-dual solution to the reduced problem,
it recovers an (approximate) primal-dual solution to the original problem.

All the explorers and reductions of cuPSLP are implemented on the GPU, while
postsolve runs on the CPU. During presolve, each reduction is logged onto a
postsolve stack that encodes the information needed to reverse it. This stack is
the one data structure whose size we cannot bound \emph{a priori}, and on
problems with many reductions it can exceed the device memory we dedicate to it.
We therefore maintain the stack on the GPU but flush it to the host as it
approaches its capacity.

\section{Experiments}
\label{sec:experiments}
To study the impact of GPU acceleration on presolving, we compare the presolve
time of cuPSLP against that of its CPU counterpart PSLP. The two presolvers
implement the same set of explorers, but in cuPSLP we have redesigned them to
better exploit GPU parallelism, so reductions are invoked in a different order.
We therefore also verify that this redesign has not hurt reduction quality. We
do not study the impact of presolve on total solve time of a GPU-accelerated
solver, but refer to \cite{Cederberg2026} for experiments showing that cuPDLPx
greatly benefits from PSLP. We also do not include a commercial presolver in our
comparison, but note that PSLP was shown in \cite{Cederberg2026} to be several
times faster than Gurobi's presolver while retaining similar reduction quality.
Comparing cuPSLP against PSLP therefore indirectly benchmarks it against a
state-of-the-art commercial presolver as well.

\subsection{Experimental setup}

\paragraph{Benchmark sets.}
We conduct numerical experiments on two LP benchmark sets. First, we
consider the Mittelmann LP benchmark set \cite{MittelmannLP}, which consists of
49 publicly available instances. Second, we consider the GAMS large-scale LP
benchmark set \cite{GamsLp2026}, which consists of 21 instances. These problems
are generally larger than those in the Mittelmann LP benchmark set.

\paragraph{Hardware.}
We benchmark the impact of GPU acceleration using a cloud virtual machine with 6
cores (12 vCPUs) of an Intel Xeon Gold 6338 CPU paired with an NVIDIA A100 (80
GB) GPU. We do acknowledge that results may vary on different hardware, but the
A100 is widely available and this node is a standard configuration offered by
GPU providers. We run PSLP on the CPU and cuPSLP on the GPU, built with CUDA
12.9.

\paragraph{Metrics.}
We measure the shifted geometric mean presolve time (with shift = 1 second) and
the ratio of the number of nonzeros in the reduced problem to that in the
original problem. The former measures speed, and the latter serves as a proxy for the
quality of the reduction. We also report
results on individual instances. We do not include the time to transfer the
original and reduced problems between host and device, but note that cuPSLP has
the advantage that the reduced problem resides on the GPU after presolve and can
be passed directly to the next GPU-resident stage of the solver pipeline without
an extra transfer.

\subsection{Numerical results}
Table \ref{tab:presolve-times} reports the shifted geometric mean presolve
times. On the Mittelmann and GAMS benchmark sets, cuPSLP on the A100 is
$11\times$ and $42\times$ faster than PSLP on the Xeon CPU, respectively. These
results corroborate our central finding that GPU acceleration can greatly speed
up presolving.

\begin{figure}[ht]
\centering
\includegraphics[width=0.85\textwidth]{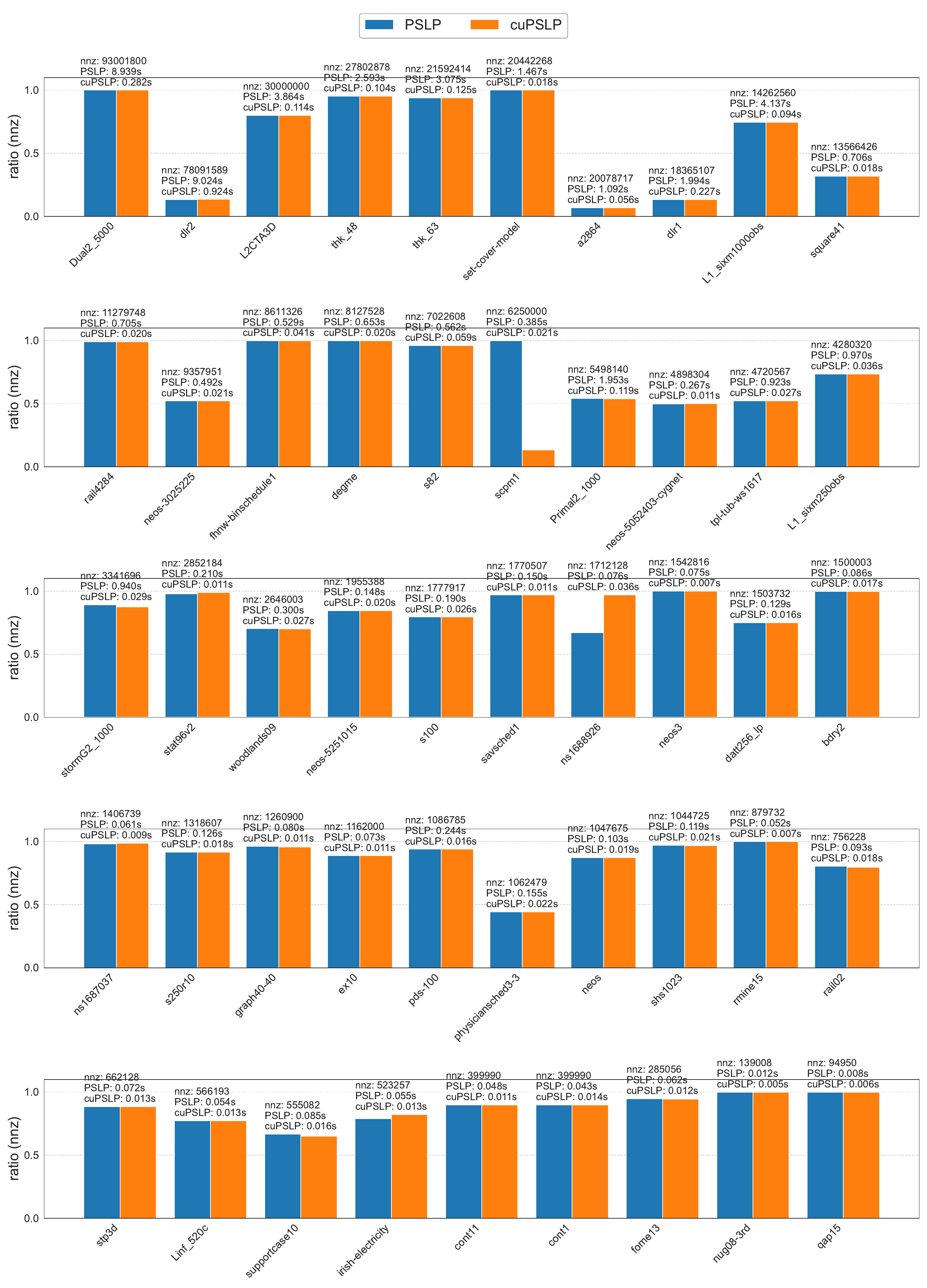}
\caption{Presolve results on the Mittelmann LP benchmark set for cuPSLP on the
A100 versus PSLP on the Xeon CPU.}
\label{fig:presolve-time-mittelmann}
\end{figure}

\begin{figure}[ht]
\centering
\includegraphics[width=0.85\textwidth]{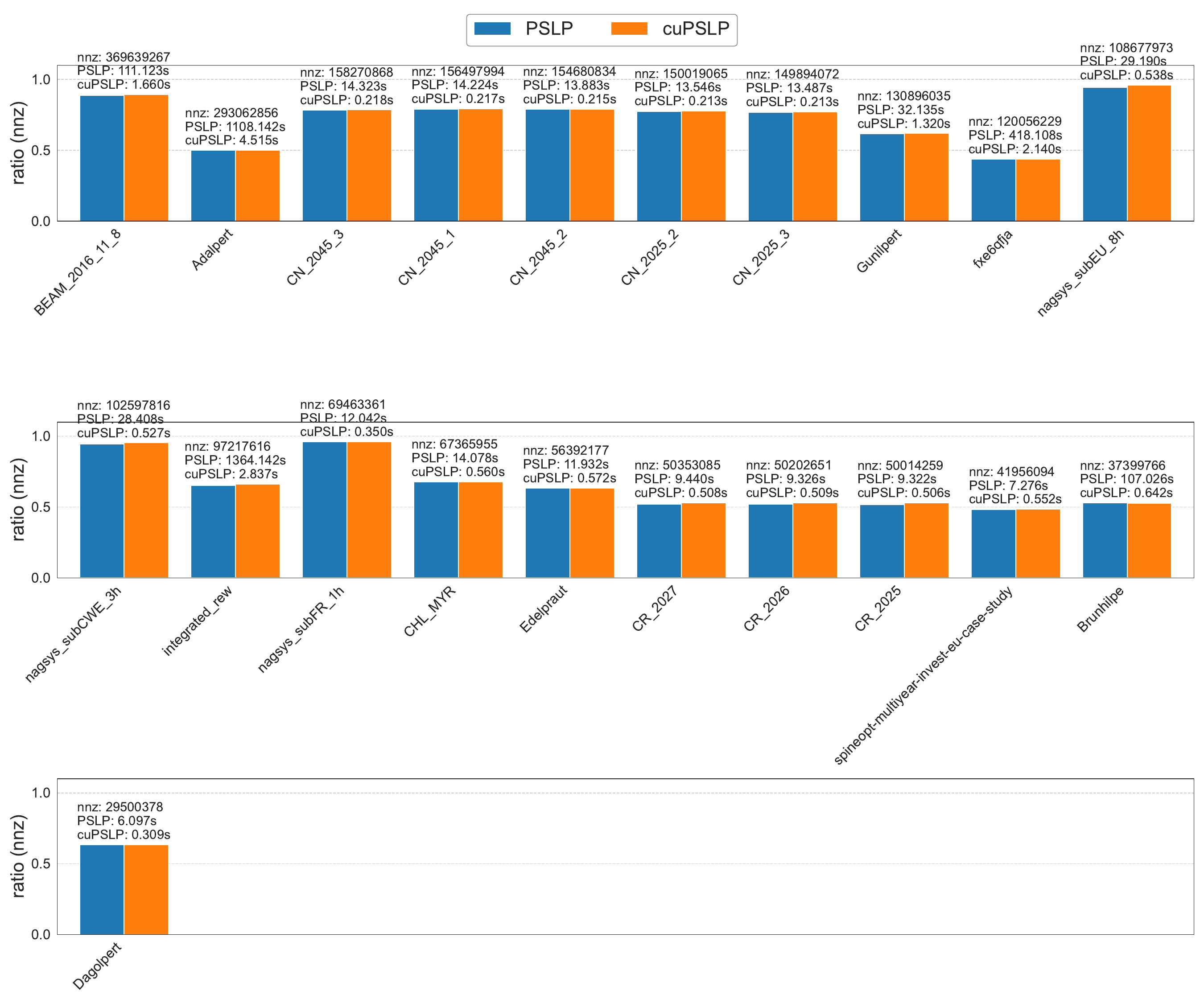}
\caption{Presolve results on the GAMS large-scale LP benchmark set for cuPSLP 
on the A100 versus PSLP on the Xeon CPU.}
\label{fig:presolve-time-gams}
\end{figure}

Figures \ref{fig:presolve-time-mittelmann} and \ref{fig:presolve-time-gams}
show, for each instance, the presolve time and the number of nonzeros in the
reduced problem for both presolvers. The reduced problems are of similar size,
except for \texttt{scpm1} where cuPSLP reduces further, and \texttt{ns1688926}
where PSLP does. We therefore conclude that the redesign of the explorers to
better exploit GPU parallelism has not hurt reduction quality in any systematic
way.

\begin{table}[t]
\centering
\begin{tabular}{lcc}
\toprule
& PSLP (CPU) & cuPSLP (GPU) \\
Dataset & Xeon & A100 \\
\midrule
Mittelmann  & 0.581 & 0.051 \\
GAMS        & 30.4 & 0.73 \\
\bottomrule
\end{tabular}
\caption{Shifted geometric mean presolve times in seconds (shift $= 1$ s).}
\label{tab:presolve-times}
\end{table}

\clearpage

\section{Discussion}
\label{sec:discussion}
We have shown that GPU acceleration can speed up presolving by more than an
order of magnitude, challenging the longstanding view that presolving is
inherently sequential and limited to the CPU. While it is exciting that another
component of a complete solver pipeline can benefit substantially from GPU
acceleration in terms of speed, we also highlight one limitation of GPU
presolve.

At least one common explorer seems inherently difficult to parallelize: the
detection and removal of linearly dependent equality constraints. Measured by
the nonzero count of the reduced problem, PSLP (and hence cuPSLP) already
recovers, without this explorer, most of the reduction achieved by a commercial
presolver \cite{Cederberg2026}, so it would likely add little on that metric.
Nonzero count, however, may not reflect the full value of an explorer. For
example, interior-point solvers that solve the KKT system via the normal
equations, whose coefficient matrix is $A D A^T$ (where $D$ is a positive
definite diagonal matrix), may be less robust when $A$ does not have full row
rank, so removing linearly dependent equality constraints can improve the
robustness of the solver even when it barely changes the size of the problem.
The classical way of detecting the dependencies, however, relies on
pivoting-based Gaussian elimination, carried out through simplex basis machinery
\cite{Andersen1995b} or via a rank-revealing sparse LU factorization of the
equality constraints \cite[\S 5.1]{Achterberg2019}, in which each pivot is
selected numerically from values produced by the eliminations before it. This is
the same sequential, data-dependent structure that has long been a challenge for
the parallelization of the simplex method \cite{Hall2010}. A natural presolve
design is therefore hybrid: move the relevant equality constraints to the CPU,
detect linear dependencies there, and return the information to the GPU for the
rest of presolve.

\section*{Acknowledgments}
We thank Rajesh Gandham, Chris Maes, and Burcin Bozkaya on the NVIDIA cuOpt team
for their feedback and for the GPU resources that made this work possible. We
also thank Zedong Peng, developer of cuPDLPx, for very helpful feedback.

\newpage
\bibliographystyle{alpha}
\bibliography{references.bib}
\end{document}